\documentclass[a4paper,12pt]{amsart}
\usepackage{amssymb}
\usepackage{ifthen}
\usepackage{graphicx}
\usepackage{float}
\usepackage{caption}
\usepackage{subcaption}
\usepackage{cite}
\usepackage{amsfonts}
\usepackage{amscd}
\usepackage{amsxtra}
\usepackage{color}

\newtheorem{thm}{Theorem}
\newtheorem{cor}{Corollary}
\newtheorem{lem}{Lemma}

\newtheorem{rem}{Remark}

\newtheorem{conj}{Conjecture}
\theoremstyle{definition}

\newtheorem{defn}{Definition}
\newtheorem{example}{Example}
\newtheorem{prob}{Problem}

\newcounter {own}
\def\theown {\thesection       .\arabic{own}}

\newenvironment{pf}[1][]{%
 \vskip 3mm
 \noindent
 \ifthenelse{\equal{#1}{}}%
  {{\slshape Proof. }}%
  {{\slshape #1.} }%
 }%
{\qed\bigskip}

\newcounter{alphabet}
\newcounter{tmp}

\newcommand{\IN}{{\mathbb N}}
\newcommand{\IC}{{\mathbb C}}
\newcommand{\ID}{{\mathbb D}}

\def\be{\begin{equation}}
\def\ee{\end{equation}}

\newcommand{\bee}{\begin{enumerate}}
\newcommand{\eee}{\end{enumerate}}

\newcommand{\blem}{\begin{lem}}
\newcommand{\elem}{\end{lem}}
\newcommand{\bthm}{\begin{thm}}
\newcommand{\ethm}{\end{thm}}
\newcommand{\bcor}{\begin{cor}}
\newcommand{\ecor}{\end{cor}}
\newcommand{\beg}{\begin{example}}
\newcommand{\eeg}{\end{example}}
\newcommand{\begs}{\begin{examples}}
\newcommand{\eegs}{\end{examples}}
\newcommand{\bdefe}{\begin{defn}}
\newcommand{\edefe}{\end{defn}}
\newcommand{\bprob}{\begin{prob}}
\newcommand{\eprob}{\end{prob}}
\newcommand{\bei}{\begin{itemize}}
\newcommand{\eei}{\end{itemize}}

\newcommand{\bcon}{\begin{conj}}
\newcommand{\econ}{\end{conj}}
\newcommand{\bcons}{\begin{conjs}}
\newcommand{\econs}{\end{conjs}}
\newcommand{\bprop}{\begin{propo}}
\newcommand{\eprop}{\end{propo}}
\newcommand{\br}{\begin{rem}}
\newcommand{\er}{\end{rem}}
\newcommand{\brs}{\begin{rems}}
\newcommand{\ers}{\end{rems}}
\newcommand{\bo}{\begin{obser}}
\newcommand{\eo}{\end{obser}}
\newcommand{\bos}{\begin{obsers}}
\newcommand{\eos}{\end{obsers}}
\newcommand{\bpf}{\begin{pf}}
\newcommand{\epf}{\end{pf}}
\newcommand{\ba}{\begin{array}}
\newcommand{\ea}{\end{array}}
\newcommand{\beq}{\begin{eqnarray}}
\newcommand{\beqq}{\begin{eqnarray*}}
\newcommand{\eeq}{\end{eqnarray}}
\newcommand{\eeqq}{\end{eqnarray*}}

\newcommand{\ds}{\displaystyle}

\def\cc{\setcounter{equation}{0}   
\setcounter{figure}{0}\setcounter{table}{0}}

\begin{document}
\bibliographystyle{amsplain}

\title[Differential Inequalities and Univalent Functions]{Differential Inequalities and Univalent Functions}

\author{Rosihan M. Ali}
\address{R. M. Ali, School of Mathematical Sciences, Universiti Sains Malaysia, 11800 USM Penang, Malaysia}
\email{rosihan@usm.my}

\author{Milutin Obradovi\'{c}}
\address{M. Obradovi\'{c}, Department of Mathematics,
Faculty of Civil Engineering, University of Belgrade,
Bulevar Kralja Aleksandra 73, 11000
Belgrade, Serbia. } \email{obrad@grf.bg.ac.rs}

\author{Saminathan Ponnusamy
}
\address{S. Ponnusamy, Department of Mathematics,
Indian Institute of Technology Madras,
Chennai-600 036, India. }

\email{samy@iitm.ac.in}

%

\begin{abstract}
Let ${\mathcal M}$ be the class of analytic functions in the unit disk $\ID$
with the normalization $f(0)=f'(0)-1=0$, and satisfying the condition
$$\left |z^2\left (\frac{z}{f(z)}\right )''+ f'(z)\left(\frac{z}{f(z)} \right)^{2}-1\right |\leq 1,
\quad z\in \ID.
$$
Functions in $\mathcal{M}$ are known to be univalent in $\ID$. In this paper, it is shown that the harmonic mean of two functions in ${\mathcal M}$ are closed, that is, it belongs again to ${\mathcal M}$. This result also holds for other related classes of normalized univalent functions.  A number of new examples of functions in $\mathcal{M}$ are shown to be starlike in $\ID$. However we conjecture that functions in $\mathcal{M}$ are not necessarily starlike, as apparently supported by other examples.
\end{abstract}


\maketitle
\pagestyle{myheadings}
\markboth{R. M. Ali, M. Obradovi\'{c}, and S.Ponnusamy}{Differential inequalities and univalent functions}
\cc

\section{Introduction }

Let ${\mathcal H}$ denote the family of analytic functions in  the open unit
disk $\ID := \{z\in \IC:\, |z| <1 \},$ and ${\mathcal A}$ its subclass of normalized functions $f(z)=z+a_2z^2+a_3z^3+\cdots$.
Further, let ${\mathcal S}$ denote the subclass of ${\mathcal A}$ consisting of functions $f$ univalent in $\ID$. Denote by $\mathcal{S}^\ast$ and $\mathcal{C}$ respectively the subclasses of $\mathcal{S}$ consisting of starlike and convex functions. Functions $f\in \mathcal{S}^*$
map $\ID$ onto starlike domains with  respect to the origin, while $f \in \mathcal{C}$ whenever $f(\ID)$ is a convex domain. Analytically, $f \in \mathcal{S}^\ast$ if ${\rm Re}\left(zf'(z)/f(z)\right) > 0$, while $f \in \mathcal{C}$ if ${\rm Re} \left( 1+ zf''(z)/f'(z)\right) > 0$.

Investigations into particular subclasses of ${\mathcal A}$ continued to be of recent interest. These include the class $\mathcal{U}$ consisting of functions $f\in \mathcal{A}$ satisfying
$$ \left| f'(z)\left(\frac{z}{f(z)} \right)^{2}-1\right| \leq 1, \quad z\in \ID,$$
as well as the class $\mathcal{P}$ of functions $f\in \mathcal{A}$ with
$$ \left|\left(\frac{z}{f(z)}\right)''\right|\leq 2, \quad z\in \ID. $$
%
The strict inclusion $\mathcal{P}\subsetneq \mathcal{U}\subsetneq
\mathcal{S}$ holds within these classes (see \cite{Ak,NOO-89,ON-72} for a proof).  There are several generalizations \cite{OP-RRMPA-09} of this result. For recent investigations on $\mathcal U$ and its generalization,
we refer to \cite{OPW,PW-2017,PW-2017b} and the references therein.

In this paper, the phrase $f\in {\mathcal{U}}$ (respectively, $ f \in {\mathcal{P}}$) in $|z|<r$ means that
the defining inequality holds in $|z|<r$ instead of the full disk $|z|<1$. We also follow this standard convention for other classes. In \cite{OP-Kodai2011} and \cite{OP-AML-12}), the authors discussed the classes ${\mathcal M}$
and ${\mathcal N}$ of functions from ${\mathcal A}$ satisfying respectively the differential inequality
$$\left |{\mathcal M}_{f}(z)\right |\leq 1, \quad \mbox{and}~ \left |{\mathcal N}_{f}(z)\right |\leq 1, \quad z\in \ID,$$
where
$${\mathcal M}_{f}(z)=z^2\left (\frac{z}{f(z)}\right )''+ f'(z)\left(\frac{z}{f(z)} \right)^{2}-1
$$
and
$$ {\mathcal N}_{f}(z)= -z^3\left (\frac{z}{f(z)}\right )'''+ f'(z)\left(\frac{z}{f(z)} \right)^{2}-1.
$$
These classes are also closely related to the class ${\mathcal U}$ in the sense of the strict inclusions
$$\mathcal{N}\subsetneq \mathcal{M}\subsetneq \mathcal{P} \subsetneq \mathcal{U}.
$$
A slightly general version of this result is given in \cite{BP-2013}.

In  \cite{OP-BMSS2012}, Obradovi\'{c}, and Ponnusamy discussed ``harmonic mean'' of two univalent analytic
functions. These are functions $F$ of the form
\be\label{5-10eq1}
F(z)=\frac{2f(z)g(z)}{f(z)+g(z)},
\ee
or equivalently,
\be\label{eq2a}
\frac{1}{F(z)}-\frac{1}{z}=\frac{1}{2}\left [\left(\frac{1}{f(z)}-\frac{1}{z}\right ) +\left (\frac{1}{g(z)}-\frac{1}{z}\right)
\right ],
\ee
where $f,g\in \mathcal{S}$. In particular,
the authors in \cite{OP-BMSS2012} determined the radius of univalency of $F$, and proposed the following two conjectures.

\bcon\label{con-new1}
\begin{enumerate}
\item[{\rm \textbf{(a)}}]
The function $F$ defined by {\rm (\ref{5-10eq1})} is not necessarily univalent in $\ID$ whenever
$f,g\in \mathcal{S}$ such that  $((f(z)+g(z))/z) \neq 0$ in  $\ID$.
\item[{\rm \textbf{(b)}}] The function $F$ defined by {\rm (\ref{5-10eq1})} is  univalent in $\ID$ whenever
$f,g\in \mathcal{C}$ such that  $((f(z)+g(z))/z) \neq 0$ in $\ID$.
\end{enumerate}
\econ
The authors in \cite{OP-BMSS2012} showed that whenever $f,g\in \mathcal{U}$, then the function $F$  defined by \eqref{5-10eq1} belongs
to $\mathcal{U}$ in the disk $|z|<\sqrt\frac{\sqrt{5}-1}{2}\approx 0.78615$.

While Conjecture \ref{con-new1} remains open, the aim of this paper is to show that Conjecture \ref{con-new1}\, \textbf{(a)} does not hold when the class $\mathcal S$ is replaced by $\mathcal U$. Indeed, it does not hold true even for the classes ${\mathcal M}, \,{\mathcal N},\, \mbox{and}~ {\mathcal P}$. The second objective of the paper is to consider several examples in examining starlikeness of functions in the classes ${\mathcal M}, \,{\mathcal N},\, \mbox{and}~ {\mathcal P}$. We conclude with a conjecture that functions in the class ${\mathcal M}$ are not necessarily starlike in $\ID$.

\section{On the harmonic mean of univalent functions}

\bthm\label{th1}
Let $f,g \in \mathcal{U}$ satisfy  $\frac{f(z)+g(z)}{z}\neq 0$ for $z\in \ID$. Then the function $F$ given by
\eqref{5-10eq1} also belongs to the class $\mathcal{U}$.
\ethm\bpf
From \eqref{eq2a}, it readily follows from the triangle inequality that the function $F$ satisfies
\beqq
\left|F'(z)\left(\frac{z}{F(z)}\right)^{2}-1\right|&=&
\left|-z^{2}\left(\frac{1}{F(z)}-\frac{1}{z}\right)'\right|
\\
&\leq & \frac{1}{2}\left|-z^{2}\left(\frac{1}{f(z)}-\frac{1}{z}\right)'\right|+\frac{1}{2}\left|
-z^{2}\left(\frac{1}{g(z)}-\frac{1}{z}\right)'\right|\\
&=&\frac{1}{2}\left|f'(z)\left(\frac{z}{f(z)}\right)^{2}-1\right|
+\frac{1}{2}\left|g'(z)\left(\frac{z}{g(z)}\right)^{2}-1\right|<1.
\eeqq
Thus $F \in \mathcal{U}$.
\epf

Moreover, we see that Theorem \ref{th1} holds true if the class $\mathcal{U}$ is replaced by
the class ${\mathcal M}$.

\bthm\label{th2}
Suppose $f,g \in \mathcal{M}$ satisfy  $\frac{f(z)+g(z)}{z}\neq 0$ for $z\in \ID$. Then the function $F$ given by
\eqref{5-10eq1} also belongs to the class $\mathcal{M}$.
\ethm\bpf
Now
$$f'(z)\left(\frac{z}{f(z)}\right)^{2}-1= -z^{2}\left(\frac{1}{f(z)}-\frac{1}{z}\right)'.
$$
Using this equality, it follows that
\beqq
{\mathcal M}_{f}(z)&=&z^{2}\left[\left (\frac{z}{f(z)}\right )''-\left(\frac{1}{f(z)}-\frac{1}{z}\right)'\right]\\
&=&z^{2}\left[\left(\left (\frac{z}{f(z)}\right )'-\frac{1}{f(z)}+\frac{1}{z}\right)'\right]\\
&=&z^{2}\left[\left(z\left (\frac{1}{f(z)}\right )'+\frac{1}{z}\right)'\right]\\
&= & z^{2}\left[\left(z\left (\frac{1}{f(z)}-\frac{1}{z}\right)'\right)'\right]\\
&=&z^{3}\left (\frac{1}{f(z)}-\frac{1}{z}\right)''+z^{2}\left (\frac{1}{f(z)}-\frac{1}{z}\right)'.
\eeqq
In view of \eqref{eq2a}, this means that
$${\mathcal M}_{F}(z)=\frac{1}{2}\left({\mathcal M}_{f}(z)+{\mathcal M}_{g}(z)\right),
$$
and use of the triangle inequality yields the desired result.
\epf

\bthm\label{th3}
Let $f,g \in \mathcal{N}$ satisfy  $\frac{f(z)+g(z)}{z}\neq 0$ for $z\in \ID$. Then the function $F$ given by
\eqref{5-10eq1} also belongs to the class $\mathcal{N}$.
\ethm\bpf
As in the proof of Theorem \ref{th2}, we see that
\beqq
{\mathcal N}_{f}(z)&=&-z^{2}\left[ z\left(\left(\frac{z}{f(z)}\right)'\right )''+\left(\frac{1}{f(z)}-\frac{1}{z}\right)'\right]\\
&=&-z^{2}\left[ z\left(\frac{1}{f(z)}
-\frac{1}{z}f'(z)\left(\frac{z}{f(z)}\right)^{2}\right )''+\left(\frac{1}{f(z)}-\frac{1}{z}\right)'\right]\\
&=&-z^{2}\left[ z\left(z \left(\frac{1}{f(z)}-\frac{1}{z}\right)'+\frac{1}{f(z)}-\frac{1}{z}\right)''
+\left(\frac{1}{f(z)}-\frac{1}{z}\right)'\right]\\
&=&-z^{4}\left(\frac{1}{f(z)}-\frac{1}{z}\right)'''-3z^{3}\left(\frac{1}{f(z)}-\frac{1}{z}\right)''
-z^{2}\left(\frac{1}{f(z)}-\frac{1}{z}\right)'.
\eeqq
Thus relation \eqref{eq2a} gives
$${\mathcal N}_{F}(z)=\frac{1}{2}\left({\mathcal N}_{f}(z)+{\mathcal N}_{g}(z)\right),
$$
and the proof of theorem readily follows.
\epf

Finally, it is also readily shown that the above theorem holds true for the class $\mathcal{P}$.

\section{Examples and a Conjecture}

It is known that functions in the class ${\mathcal U}$ are not necessarily starlike. There are a number of examples displaying functions in ${\mathcal U}$ that are not starlike in $\ID$, see for instance \cite{OP-01}. However, is ${\mathcal M} \subset {\mathcal S}^*$? This section discusses the latter problem.
{\rm
\beg\label{KMJ-eg1}
To present a one-parameter family of functions in ${\mathcal M}$ that are also starlike, consider the function $f$ given by
$$ \frac{z}{f(z)}=1+(1-\alpha)z+\alpha z^m,
$$
where $\alpha\in (0,1)$ and $m\in\IN\backslash \{1\}=\{2,3, \ldots \}$ are such that $\alpha (m-1)^2=1.$
Then $z/f(z)\neq 0$ in $\mathbb{D}$ and
$$\sum_{k=2}^{\infty}(k-1)^2|b_k|=(m-1)^2\alpha=1,
$$
and therefore, $f \in {\mathcal M}.$

Next, we show that $f$ is starlike whenever $m>1$ is an odd integer. Now,
a simple calculation shows
$$\frac{zf'(z)}{f(z)}=\frac{1-\alpha(m-1)z^m}{1+(1-\alpha)z+\alpha z^m}.
$$
With $z=e^{i\theta},$ then
$$\frac{e^{i\theta}f'(e^{i\theta})}{f(e^{i\theta})}=\frac{A(\theta)+iB(\theta)}{|1+(1-\alpha)e^{i\theta}+\alpha e^{im\theta}|^2},
$$
where
\beqq
A(\theta) &= &1+(1-\alpha)\cos\theta-\alpha(m-2)\cos(m\theta)
\\
& &~~~~ -\, \alpha(1-\alpha)(m-1)\cos(m-1)\theta-\alpha^2(m-1).
\eeqq
Note that $A(\theta)=A(-\theta)$.
As $ \alpha= 1/(m-1)^2,$ the expression for $A(\theta )$ reduces to
$$ A(\theta)=1-\frac{1}{(m-1)^3}-\frac{m(m-2)}{(m-1)^2} D(\theta ),
$$
where
$$ D(\theta)=-\cos\theta +\frac{1}{m}\cos(m\theta)+\frac{\cos(m-1)\theta}{m-1}.
$$

To show starlikeness, that is, $f\in \mathcal{S}^*$, it suffices to show that $A(\theta)\geq0$ for
$0\leq\theta\leq\pi.$
First we prove the assertion for the case $m=3$, while the
general case is obtained separately. Setting $m=3$, $A(\theta)$ reduces to
$$A(\theta)= \frac{7}{8}-\frac{3}{4}\left [-\cos\theta +\frac{1}{3}\cos  3\theta +\frac{1}{2}\cos  2\theta \right ],
$$
and from the identities $\cos 2\theta =2\cos ^2\theta -1$ and $\cos 3\theta =4\cos ^3\theta -3\cos \theta $,
\beqq
A(\theta)&=& \frac{1}{4}(5+6\cos \theta -4\cos ^3 \theta -3 \cos^2 \theta)\\
&=&\frac{1}{4}(1+\cos \theta)^2(5-4\cos \theta),
\eeqq
which shows  that $A(\theta)\geq 0$. Thus, the function $f_3(z)$ given by
$$f_3(z)=\frac{z}{1+\frac{3}{4} z+\frac{1}{4} z^{3}}
=\frac{4z}{(1+z)(4-z+z^2)},$$
is starlike in $\mathbb D$.

Next, we proceed to prove starlikeness for the general case. This requires more computations. First,
\begin{align*}
D'(\theta)&=\sin\theta-\sin(m\theta)-\sin(m-1)\theta\\
&= \sin\theta-2\sin\frac{(2m-1)\theta}{2}\cos\frac{\theta}{2}\\
&=2\cos\frac{\theta}{2}\left[\sin\frac{\theta}{2} {\color{red}{-}}\sin\frac{(2m-1)\theta}{2}\right]\\
&= 4\cos\frac{\theta}{2}\cos\frac{m\theta}{2}\sin\frac{(m-1)\theta}{2}.
\end{align*}
We need to show that $A(\theta)\geq 0$ for $0\leq\theta\leq\pi$.
It is convenient to set $m=2n+1$, $n\ge 2$ so that
$$D'(\theta) =4\cos\frac{\theta}{2}\cos\frac{(2n+1)\theta}{2}\sin n\theta, \quad n\ge 2,
$$
where $D(\theta)$ takes the form
$$ D(\theta)=-\cos\theta +\frac{1}{2n+1}\cos(2n+1)\theta+\frac{1}{2n}\cos (2n\theta).
$$
Clearly, $D'(\theta)=0$ for $\theta =0, \pi$, and the critical points of $D(\theta)$ in the open interval $(0,\pi)$
are given by
$$\left \{ \ba{lll} \theta _j  =& \ds  \frac{(2j-1)\pi}{2n+1} &\mbox{for $j=1,2,\ldots, n$},\\[3mm]
 \theta _j' =& \ds \frac{j\pi}{n} & \mbox{for $j=1,2,\ldots, n-1$,}
 \ea
 \right .
$$
$n\ge 2$. Moreover, for each $n\ge 2$,
$$\left \{ \ba{rl}
\ds \cos\frac{(2n+1)\theta}{2} > 0 & \ds\mbox{for } 0< \theta <\theta _1, \\[2mm]
\ds (-1)^{j}\cos\frac{(2n+1)\theta}{2} >0 &\mbox{for $ \theta _j < \theta <\theta _{j+1}$ and  for  $j=1,2,\ldots, n$},\\[3mm]
\ds (-1)^{j-1}\sin n\theta >0 & \ds \mbox{for $ \theta _{j-1}' < \theta <\theta _{j}'$ and for $j=1,2,\ldots, n$}.
 \ea
 \right .
$$

In view of the above inequalities and after a careful scrutiny, it follows that
$$ D'(\theta) \left \{ \ba{rl}
= 0  &\mbox{for $\ds  \theta =0, \theta _j, \theta _{j}'$ for $j=1,2,\ldots, n$},\\
>0 &\mbox{for $\ds  \theta \in (0 , \theta _{1})\cup(\theta _j' , \theta _{j+1})$  for $j=1,2,\ldots, n-1$},\\
<0 &\mbox{for $\ds  \theta \in (\theta _j , \theta _{j}')$  for $j=1,2,\ldots, n$},
\ea
\right .
$$
where $0<\theta _{1}<\theta _{1}'<\theta _{2}<\cdots <\theta _{j}<\theta _{j}'<\theta _{j+1}<\cdots <\theta _{n}<\theta _{n}'=\pi$.
Therefore,
$$ D(\theta)   \leq  \max\left \{ D(0), D\left (\theta _j \right ), D\left (\theta _j' \right ):\,j=1,2,\ldots , n\right \}.
$$
Since
\begin{align*}
D(0)& =-1+\frac{1}{2n+1} + \frac{1}{2n}=-\frac{2n}{2n+1} + \frac{1}{2n},\\
D(\pi)& =1-\frac{1}{2n+1} + \frac{1}{2n} =\frac{2n}{2n+1} + \frac{1}{2n}>0,
\end{align*}
then $D(0)\leq D(\pi)$. Moreover,
\begin{align*}
D\left (\theta _j \right ) &=-\cos \theta _j +\frac{1}{2n+1} \cos (2j-1)\pi
+\frac{1}{2n} \cos (2n+1-1)\theta _j\\
&=-\cos \theta _j -\frac{1}{2n+1}- \frac{1}{2n} \cos \theta _j\\
&=-\left (\frac{2n+1}{2n} \right )\cos \theta _j -\frac{1}{2n+1},
\end{align*}
and
\begin{align*}
D\left (\theta _j' \right ) &=-\cos \theta _j' +\frac{1}{2n+1} \cos (2n+1)\frac{j}{n}\pi
+\frac{1}{2n} \cos (2j\pi)\\
&=-\left (1-\frac{1}{2n+1} \right )\cos \theta _j' +\frac{1}{2n}\\
&=-\frac{2n}{2n+1}\cos \theta _j' +\frac{1}{2n}.
\end{align*}
We deduce that $D\left (\theta _j \right )\leq D(\pi) $ and $D\left (\theta _j' \right )\leq D(\pi)$ holds
for each  $j=1,2,\ldots , n$. Thus,
$D(\theta )\leq D(\pi)$ for $\theta \in [0,\pi]$. This observation shows that
$$A(\theta)\geq A(\pi) = 1-\frac{1}{8n^3} -\frac{(2n+1)(2n-1)}{4n^2}\left (\frac{2n}{2n+1} +\frac{1}{2n}
\right )
=0 ~\mbox{  for $\theta \in [0,\pi]$}.
$$
Hence ${\rm Re}\,(e^{i\theta}f'(e^{i\theta})/f(e^{i\theta}))\geq 0$,  which implies that $f$ is starlike in $\ID.$
Summarizing, for each $n\geq 1$, the function $f_n$ given by
$$\frac{z}{f_n(z)}=1+\left (1-\frac{1}{4n^2}\right )z+\frac{1}{4n^2} z^{2n+1},
$$
belongs ${\mathcal M}$, and $f_n$ is starlike in $\ID$.
\eeg
}
\beg
{\rm
Consider
$$ f(z)=\frac{z}{\phi (z)}, \quad \phi (z)=1+\left (1-\frac{\zeta(5)}{\zeta(3)}\right )z
+\frac{1}{\zeta(3)}\sum_{n=2}^{\infty}\frac{z^{n}}{(n-1)^5}.
$$
We may rewrite $\phi$ as
$$\phi (z)=
 1+\left (1-\frac{\zeta(5)}{\zeta(3)}\right )z+\frac{1}{\zeta(3)}\frac{z^2}{4!}\int_{0}^{1}
\frac{(\log(1/t))^4\,dt}{1-tz}.
$$
It is a simple exercise to see that $\phi (z)\neq 0$ in $\ID$ and $f\in {\mathcal M}$. The Mathematica software is used to display the
image of the unit disk under $f$ as shown in Figure \ref{fig1}. It apparently displays that $f(\ID)$ is a starlike domain.
\begin{figure}
\includegraphics[scale=1.4]{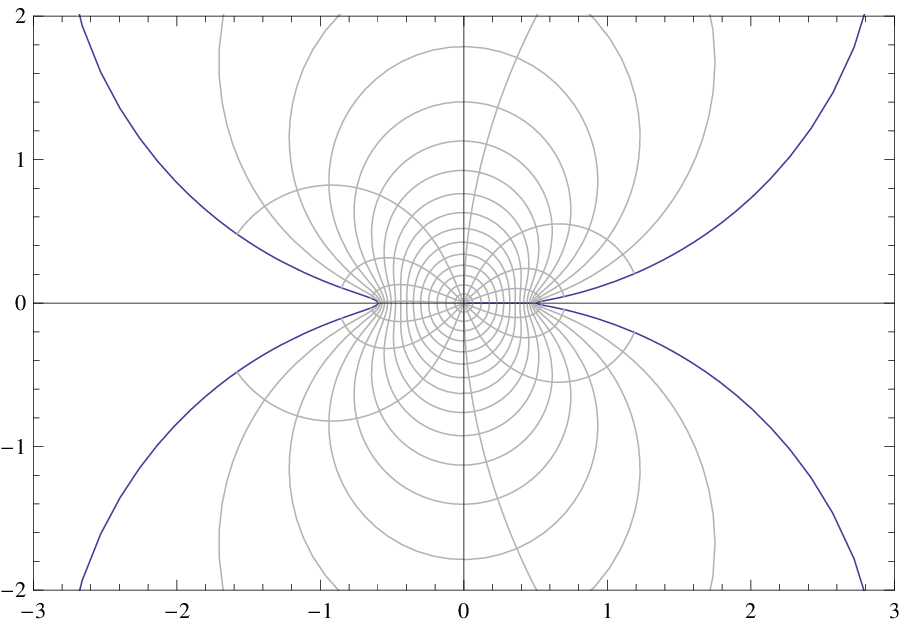}
\caption {\label{fig1}}
\end{figure}
}
\eeg

\beg
{\rm It is illustrative to present a general example showing that functions in ${\mathcal U}$ do not necessarily belong to
${\mathcal S}^*$. For $n\geq 3$, consider the function
$$f_n(z)=\frac{z}{1+ibz+(1/(n-1))e^{2i\beta}z^n}.
$$
For $|b|\leq (n-2)/(n-1)$ and $\beta$ a real number, then
$$ {\rm Re}\, \left (\frac{z}{f_n(z)}\right ) > 1-|b|-\frac{1}{n-1} \geq 0,
$$
and
$$
\left |\left (\frac{z}{f_n(z)} \right ) ^2f_n'(z)-1\right |
= \left |-e^{2i\beta}z^n \right | <1~\mbox{ for $z\in \ID$},
$$
so that $f_n\in {\mathcal U}$ for each $n\geq 3$. On the other hand, $f_n$ is not in ${\mathcal S}^*$ when
$0<b\leq (n-2)/(n-1)$ and $0<\beta<\arctan (b(n-1)/(n-2))$. This follows on account that
$$ \left .{\rm Re}\, \left (\frac{zf_n'(z)}{f_n(z)}\right ) \right |_{z=1} =
 \frac{[ (2(n-2)/(n-1))\sin \beta  -2b\cos \beta ]\sin \beta }{| 1+ib+(e^{2i\beta}/(n-1))|^2}< 0.
$$
}
\eeg

\beg\label{KMJ-eg4}
{\rm
Consider the function $f$ defined by
$$ \frac{z}{f(z)}=1+(1-\alpha)z+\alpha z^{m},
$$
where $\alpha\in (0,1)$ and $m\geq 3$ is an odd integer such that $\alpha m(m-1)=2.$
Then $z/f(z)\neq 0$ in $\mathbb{D}$ and
$$ \left |\left (\frac{z}{f(z)} \right )''\right |=|\alpha m(m-1)z^{m-2}|<\alpha m(m-1)=2,
$$
and therefore, $f \in {\mathcal P}.$

As in Example \ref{KMJ-eg1},
$${\rm Re}\left (\frac{e^{i\theta}f'(e^{i\theta})}{f(e^{i\theta})}\right )
=\frac{A(\theta )}{|1+(1-\alpha)e^{i\theta}+\alpha e^{im\theta}|^2},
$$
where
\beqq
A(\theta) &= &1+(1-\alpha)\cos\theta-\alpha(m-2)\cos(m\theta)
\\
& &~~~~ -\, \alpha(1-\alpha)(m-1)\cos(m-1)\theta-\alpha^2(m-1).
\eeqq
Substituting $\alpha= 2/(m(m-1))$ and $m=2n+1$ ($n\ge 1$), the last expression for $A(\theta )$ reduces to
\be\label{new-atheta1}
A(\theta)=
1-\frac{2}{n(2n+1)^2}+\frac{2n-1}{n(2n+1)}D(\theta),
\ee
where
$$D(\theta)=(n+1)\cos\theta -\cos(2n+1)\theta -
\frac{2(n+1)}{2n+1}\cos 2n\theta.
$$

To prove that $f$ is not starlike in $\ID$, it suffices to show that $A(\theta)<0$
for some $\theta \in (-\pi,\pi)$. In the case of $m=3$ (i.e. $n=1$),  it is a simple
exercise to see that
$$A(\theta) =\frac{1}{9}(1+\cos \theta)(11+4\cos \theta -12\cos ^2\theta ),
$$
which is clearly negative for $\theta$ near $\pi$. Indeed, substituting $\cos \theta =-8/9$
or $\theta _0 =6\pi/7$, it can be verified that $A(\theta) \approx  -55/2187<0$,
and $A(\theta_0) \approx  -0.25811<0$. Thus, the function
$$f_3(z)=\frac{z}{1+\frac{2}{3}z+\frac{1}{3} z^{3}}=\frac{3z}{(1+z)(3-z+z^2)}
$$
belongs to ${\mathcal P}\backslash {\mathcal S}^*$.

To do away the problem for some other values of $n$, we proceed as follows. Set
$$\theta =\frac{2(2n+1)\pi}{4n+3}~ \mbox{ and }~\phi =\frac{\pi}{2(4n+3)}
$$
so that $\phi =(\pi-\theta)/2$.
Then
$$\cos \theta =-\cos 2\phi = 2\sin^2\phi -1,
$$
$$\cos (2n+1)\theta =-\cos 2(2n+1)\phi =-\sin \phi,
$$
and
$$\cos 2n\theta =\cos 4n\phi =\sin  3\phi =3\sin \phi -4\sin ^3\phi.
$$
Thus, $A(\theta)$ given by  (\ref{new-atheta1}) can be simplified leading to
\begin{align*}
 A(\theta) =~&
1-\frac{2}{n(2n+1)^2} -\frac{2(2n-1)(n+1)}{2n(2n+1)}\\
 &~~~~+\frac{2n-1}{n(2n+1)}\left [ 2(n+1)\sin^2\phi -\frac{4n+5}{2n+1}\sin\phi +\frac{8(n+1)}{2n+1}\sin^3\phi
\right ].
\end{align*}
It is seen from the computer algebra system Mathematica that $A (\theta)<0$ for $n=1,2, \ldots, 15$. For easy reference, Table \ref{table1} lists the values of $A(\theta)$ for $n=1,2, \ldots, 14$.
\begin{table}
\center{
\begin{tabular}{|c|c||c|c|c|c|}
\hline
$n$ & value of $A(\theta)$ &$n$ & value of $A(\theta)$\\
\hline
1& $-0.0258011$ &8 & $-0.000243709$
\\
2& $-0.0103986$ &9 & $-0.000154718$
\\
3& $-0.00437311$ &10 & $-0.0000989276$
\\
4& $-0.00211511$ &11 & $-0.0000628326$
\\
5& $-0.00113174$ &12 & $-0.0000388937$
\\
6& $-0.00064961$ &13  & $-0.000022708$
\\
7& $-0.00039145$ & 14  & $-0.0000116051$
\\
\hline
\end{tabular}
}
\bigskip
\caption{Values of $A(\theta )$ for certain choices of $\theta$ \label{table1}
}
\end{table}
Thus, we conclude that the above procedure helps us to show that for each
$n\in \{1, 2, \ldots, 14\}$, the function $f_n$ given by
$$\frac{z}{f_n(z)}=1+\left (1-\frac{1}{n(2n+1)}\right )z+\frac{1}{n(2n+1)} z^{2n+1}
$$
is not starlike in $\ID$. By a minor modification in the choice of $\theta$, one can show
that $f_n$ is not starlike for some $n\ge 15$ although it is not clear whether $f_n$ is starlike
for larger values of $n$.
}
\eeg

The ideas and the motivations behind the above examples lead to the following conjecture:

\bcon
The class $\mathcal M$ is not contained in ${\mathcal S}^*$.
\econ

\subsection*{Acknowledgments}
The first author gratefully acknowledged support from a Universiti Sains Malaysia research university grant 1001/PMATHS/8011101. The  work of the third author is supported by Mathematical Research Impact Centric Support of DST, India  (MTR/2017/000367).



\begin{thebibliography}{150}

\bibitem{BP-2013} A. Baricz, and S. Ponnusamy,
\textit{Differential Inequalities and Bessel Functions},
J. Math. Anal. and Appl. {\bf 400}(2) (2013), 558--567.

\bibitem {Ak} L.A. Aksentiev,
\textit{Sufficient conditions for univalence of regular functions. (Russian)}
\textrm{Izv. Vys\v s. U\v cebn. Zaved. Matematika}  {\bf 3}(4) (1958), 3--7.

\bibitem {FR-2006} R. Fournier and S. Ponnusamy,
\textit{A class of locally univalent functions defined by a differential inequality},
\textrm{Complex Var. Elliptic Equ.} {\bf 52}(1) (2007), 1--8.

\bibitem{Fr1}  B. Friedman,
\textit{Two theorems on schlicht functions,}
\textrm{Duke Math. J.} {\bf 13} (1946), 171--177.

\bibitem {NOO-89}  M. Nunokawa, M. Obradovi\'c, and S. Owa,
\textit{One criterion for univalency,}
\textrm{Proc. Amer. Math. Soc.} {\bf 106} (1989), 1035--1037.

\bibitem {OP-01}  M. Obradovi\'c and S. Ponnusamy,
\textit{New criteria and distortion theorems for univalent functions,}
\textrm{Complex Variables Theory Appl.}  {\bf 44} (2001), 173--191.

\bibitem{OP-RRMPA-09} M. Obradovi\'c and S. Ponnusamy,
\textit{On certain subclasses of univalent functions and radius properties},
\textrm{Rev. Roumaine Math. Pures Appl., } {\bf 54}(4) (2009), 317--329.

\bibitem{OP-Kodai2011} M. Obradovi\'c and S. Ponnusamy,
\textit{A class of univalent functions defined by a differential inequality},
Kodai Math. J.  {\bf 34} (2011), 169-–178.


\bibitem{OP-AML-12} M. Obradovi\'{c}, S. Ponnusamy,
\textit{On a class of univalent functions},
\textrm{Appl. Math. Lett.} {\bf 25}(10) (2012) 1373--1378.

\bibitem{OP-BMSS2012} M. Obradovi\'c and S. Ponnusamy,
\textit{On harmonic combination of univalent functions},
\textrm{Bull. Belg. Math. Soc. (Simon Stevin)} {\bf 19}(3) (2012),
461--472.

\bibitem{OPW} M. Obradovi\'c, S. Ponnusamy, and K.-J. Wirths,
\textit{Geometric studies on the class $\mathcal{U}(\lambda)$},
\textrm{Bull. Malaysian Math. Sci. Soc.} \textbf{39}(3) (2016), 1259--1284.




\bibitem{PW-2017}  { S. Ponnusamy and K.-J. Wirths},
\textit{Elementary considerations for classes of meromorphic univalent functions,}
\textrm{Lobachevskii J. Math.} \textbf{39}(5) (2018), 713--716.
Preprint.

\bibitem{PW-2017b}  { S. Ponnusamy and K.-J. Wirths},
 \textit{Coefficient problems  on the class $U(\lambda)$},
 \textrm{Probl. Anal. Issues Anal.} \textbf{7}(25), No. 1, (2018), 87--103.

\bibitem{ON-72}  S. Ozaki and M. Nunokawa,
\textit{The Schwarzian derivative and univalent functions,}
\textrm{Proc. Amer. Math. Soc.} {\bf 33} (1972), 392--394.

\end{thebibliography}
\end{document}